\documentclass[12pt,reqno]{amsart}

\usepackage{amsmath, amsthm, amscd, amsfonts, amssymb, graphicx, color}
\usepackage[bookmarksnumbered,plainpages]{hyperref}
\input{mathrsfs.sty}

\newtheorem{theorem}{Theorem}[section]

\newtheorem{proposition}[theorem]{Proposition}

\theoremstyle{definition}

\newtheorem{example}[theorem]{Example}

\theoremstyle{remark}
\newtheorem{remark}[theorem]{Remark}

\numberwithin{equation}{section}

\newcommand{\ip}[2]{\left\langle#1,#2\right\rangle}

\begin{document}

\title[A characterization of Hilbert $C^*$-modules]{A characterization of Hilbert $C^*$-modules over finite dimensional $C^*$-algebras}

\author[Lj. Aramba\v si\' c, D. Baki\' c, M.S. Moslehian]{Ljiljana Aramba\v si\' c $^1$, Damir Baki\' c $^2$ and Mohammad Sal Moslehian $^3$}

\address{$^{1}$ Department of Mathematics, University of Zagreb, Bijeni\v cka cesta 30, 10000 Zagreb, Croatia.}
\email{arambas@math.hr}

\address{$^{2}$ Department of Mathematics, University of Zagreb, Bijeni\v cka cesta 30, 10000 Zagreb, Croatia.}
\email{bakic@math.hr}

\address{$^3$ Department of Pure Mathematics, Ferdowsi University of Mashhad, P. O. Box
1159, Mashhad 91775, Iran; \newline Center of Excellence in Analysis
on Algebraic Structures (CEAAS), Ferdowsi University of Mashhad,
Iran.} \email{moslehian@ferdowsi.um.ac.ir and moslehian@ams.org}

\subjclass[2000]{Primary 46L08; Secondary 46L05, 46L10}
\keywords{Hilbert $C^*$-module, finite dimensional $C^*$-algebra,
$C^*$-algebra of compact operators.}

\begin{abstract}
We show that the unit ball of a full Hilbert $C^*$-module is
sequentially compact in a certain weak topology if and only if the
underlying $C^*$-algebra is finite dimensional. This provides an
answer to the question posed in J.~Chmieli\'nski et al [Perturbation
of the Wigner equation in inner product $C^*$-modules, J. Math.
Phys. 49 (2008), no. 3, 033519].
\end{abstract}

\maketitle


\section{Introduction and preliminaries}

\noindent Let ${\mathscr A}$ be a $C^*$-algebra. A linear space
${\mathscr M}$ that is an algebraic left ${\mathscr A}$-module with
$\lambda(ax)=a(\lambda x)= (\lambda a)x$ for $x\in {\mathscr M}$,
$a\in {\mathscr A}$, $\lambda\in {\mathbb C}$, is called a
\emph{pre-Hilbert ${\mathscr A}$-module} (or an \emph{inner product
${\mathscr A}$-module}) if there exists an ${\mathscr A}$-valued
inner product on ${\mathscr M}$, i.e., a mapping $\ip{\cdot}{\cdot}
\colon {\mathscr M}\times {\mathscr M}\to {\mathscr A}$ satisfying
\begin{enumerate}
\item[(i)] $\ip{\lambda x + y}{z} = \lambda \ip{x}{z}+
\ip{y}{z};$
\item[(ii)] $\ip{ax}{y} = a\ip{x}{y};$
\item[(iii)] $\ip{x}{y}^*=\ip{y}{x};$
\item[(iv)] $\ip{x}{x}\geq 0;$
\item[(v)] $\ip{x}{x} =0 \Leftrightarrow x=0,$
\end{enumerate}
for all $x, y, z \in {\mathscr M}$, $a\in {\mathscr A}$, $\lambda\in
{\mathbb C}$. Conditions (i) and (iii) yield the fact that the inner
product is conjugate-linear with respect to the second variable. It
follows from the definition that $\|x\|_{\mathscr
M}:=\sqrt{\|\ip{x}{x}\|_{\mathscr A}}$ is a norm on ${\mathscr M}$,
whence ${\mathscr M}$ becomes a normed left ${\mathscr A}$-module. A
pre-Hilbert ${\mathscr A}$-module ${\mathscr M}$ is called a
\emph{Hilbert $C^*$-module} if it is complete with respect to the
this norm. We say that a Hilbert ${\mathscr A}$-module ${\mathscr
M}$ is \emph{full} if the linear subspace $\ip{\mathscr M}{\mathscr
M}$ of ${\mathscr A}$ generated by $\{\ip{x}{y}:x,y\in {\mathscr
M}\}$ is dense in ${\mathscr A}.$ The simplest examples are usual
Hilbert spaces as Hilbert ${\mathbb C}$-modules, and $C^*$-algebras
as Hilbert $C^*$-modules over themselves via $\ip{a}{b}=ab^*.$

The concept of a Hilbert $C^*$-module has been introduced by
Kaplansky \cite{KAP} and Paschke \cite{pas}. For more information we
refer the reader e.g.\ to monographs \cite{LAN, M-T}.


Despite a formal similarity of definitions, it is well known that
Hilbert $C^*$-modules may lack many properties familiar from Hilbert
space theory. In fact, it turns out that properties of a
$C^*$-module reflect (or originate from) the properties of the
underlying $C^*$-algebra.

A particularly well behaved class is the class of Hilbert
$C^*$-modules over $C^*$-algebras of compact operators. There are
several nice characterizations of such modules (see e.g. \cite{a,
Fr, F-S, Mag, Sch}). In our proofs we make use of orthonormal bases
which exist only in Hilbert $C^*$-modules over $C^*$-algebras of
compact operators (see \cite{a,b-g1}). Recall that a system of
vectors $\{\varepsilon_i:i\in I\}$ in a Hilbert $\mathscr A$-module
${\mathscr M}$ is said to be an \emph{orthonormal basis for
${\mathscr M}$} if it satisfies the following conditions:
\begin{enumerate}
    \item $p_i:=\ip{\varepsilon_i}{\varepsilon_i}\in \mathscr A$ is a projection such that
    $p_i\mathscr A p_i={\mathbb C}p_i$ for every $i\in I;$
    \item $\ip{\varepsilon_i}{\varepsilon_j}=0$ for every $i,j\in I,i\neq j;$
    \item $\{\varepsilon_i:i\in I\}$ generates a norm-dense submodule of ${\mathscr M}.$
\end{enumerate}
If $\{\varepsilon_i:i\in I\}$ is an orthonormal basis (with the
above properties (1), (2), (3)) for ${\mathscr M}$ then the
reconstruction formula $x=\sum_{i\in I}\ip{x}{\varepsilon_i}
\varepsilon_i$ holds for every $x\in {\mathscr M},$ with the norm
convergence. Since all orthonormal bases for a Hilbert $\mathscr
A$-module ${\mathscr M}$ have the same cardinality (see
\cite{b-g1}), it makes sense to define the \emph{orthogonal
dimension} of ${\mathscr M},$ denoted by $\dim_{\mathscr A}
{\mathscr M},$ as the cardinal number of any of its orthonormal
bases.

Various specific properties of Hilbert $C^*$-modules turn out to be
particularly useful in applications. An interesting example of
investigations of this type is a recent study of the stability of
Wigner equation (see \cite{C-I-M-S} and the references therein). In
particular, the main result in \cite{C-I-M-S} is obtained for
Hilbert $C^*$-modules satisfying the following condition:

\begin{enumerate}
\item[{\sf [H]}] For each norm-bounded sequence $(x_n)$ in
${\mathscr M}$, there exist a subsequence $(x_{n_k})$ of $(x_n)$
and an element $x_0\in{\mathscr M}$ such that the sequence
$(\ip{x_{n_k}}{y})$ converges  to $\ip{x_0}{y}$ in norm for any
$y\in {\mathscr M}$.
\end{enumerate}

\vspace{.1in}

Notice that in case of a Hilbert space condition {\sf [H]} is
clearly satisfied: this is simply the fact that the unit (and hence
each) ball in a Hilbert space is weakly sequentially compact.

It is proved in \cite[Proposition 2.1]{C-I-M-S} that a Hilbert
${\mathscr A}$-module ${\mathscr M}$ satisfies condition {\sf [H]}
whenever the underlying $C^*$-algebra is finite dimensional. In this
note we prove the converse, i.e., we show that condition {\sf [H]}
is an exclusive property of the class of Hilbert $C^*$-modules over
finite dimensional $C^*$-algebras. In this way we obtain a new
characterization of such modules and answer a question posed in
\cite{C-I-M-S} concerning condition {\sf [H]}.

\vspace{.1in}

\section{The result}
For a Hilbert space $H$ we denote by $\Bbb B(H)$ and $\Bbb K(H)$
the $C^*$-algebras of all bounded, resp. compact operators acting
on $H$. We begin with a proposition that reduces the discussion to
the class of $C^*$-algebras of compact operators.

\begin{proposition}\label{p1}
Suppose that ${\mathscr M}$ is a full Hilbert $C^*$-module over a
$C^*$-algebra ${\mathscr A}$, which satisfies condition {\sf [H]}.
Then ${\mathscr A}$ is isomorphic to a $C^*$-algebra of (not
necessarily all) compact operators acting on some Hilbert space.
\end{proposition}
\begin{proof} Let us fix $y \in {\mathscr M}$.
Consider the map $T_y :{\mathscr M} \rightarrow {\mathscr M}$ given
by $T_y(x)=\langle y,x\rangle y.$ Obviously, $T_y$ is a bounded
anti-linear operator.

Let $(x_n)$ be a norm-bounded sequence in ${\mathscr M}$ and let
$(x_{n_k})$ be a subsequence of $(x_n)$ such that, for some $x_0\in
{\mathscr M},$ $\lim_{k\to\infty} \ip{x_{n_k}}{y}= \ip{x_0}{y}$ for
all $y\in {\mathscr M}.$ Then $\lim_{k \rightarrow \infty}\langle
y,x_{n_k}\rangle y =\langle y,x_0\rangle y$ for all $y \in {\mathscr
M}$. This can be restated in the following way: for each
norm-bounded sequence $(x_n)$ in ${\mathscr M}$, the sequence
$(T_y(x_n))$ has a convergent subsequence. Hence, $T_y$ is a compact
operator. Moreover, by the hypothesis, this is true for each $y \in
{\mathscr M}$.

By \cite[Proposition~1]{a}, $(4)\Rightarrow(1),$ there is a faithful
representation $\pi:{\mathscr A}\to\Bbb B(H)$ of ${\mathscr A}$ on
some Hilbert space $H$ such that $\pi(\langle y,y\rangle)\in\Bbb
K(H).$ This holds for every $y\in {\mathscr M},$ so, by
polarization,  $\pi(\langle x,y\rangle)\in\Bbb K(H)$ for all $x,y\in
{\mathscr M},$ and therefore $\pi({\mathscr A})\subseteq \Bbb K(H).$
\end{proof}

\vspace{.1in} By the preceding proposition, condition {\sf [H]} can
only be satisfied in Hilbert $C^*$-modules over $C^*$-algebras of
compact operators. (Here, and in the sequel, we identify ${\mathscr
A}$ with $\pi({\mathscr A})$, where $\pi$ is the representation from
the preceding proof.) However, even if the underlying algebra is a
$C^*$-algebra of compact operators, one still cannot conclude that
condition {\sf [H]} is satisfied.

We demonstrate this fact in the following two examples.

\begin{example}
Consider a separable infinite dimensional Hilbert space $H$ with
an orthonormal basis $(\varepsilon_n)$. We shall regard $\Bbb
K(H)$ as a Hilbert $C^*$-module over itself via the inner product
$\langle a,b\rangle=ab^*$. Let us show that $\Bbb K(H)$ does not
satisfy {\sf [H]}.

For $n \in \Bbb N$, denote by $p_n$ the orthogonal projection to
$\mbox{span}\{\varepsilon_1, \ldots , \varepsilon_n\}$. Obviously,
the sequence $(p_n)$ is norm-bounded.

Suppose that there exist a subsequence $(p_{n_k})$ and a compact
operator $a \in \Bbb K(H)$ such that $\lim_{k \rightarrow
\infty}\langle p_{n_k}, y \rangle =\langle a,y\rangle$ for all $y
\in \Bbb K(H)$. This means $p_{n_k}y^* \rightarrow ay^*$ for all
$y\in \Bbb K(H)$, which in turn gives us $p_{n_k}y\xi \rightarrow
ay\xi$ for all $y\in \Bbb K(H)$ and for all $\xi \in H$. In
particular, for every $n \in \Bbb N$, we can take $y=p_{n}-p_{n-1}$
(that is the orthogonal projection to the one-dimensional subspace
spanned by $\varepsilon_n$) and $\xi =\varepsilon_n$. Then the
preceding relation yields $a\varepsilon_n=\varepsilon_n$ for all
$n\in \Bbb N;$ i.e., $a$ is the identity operator. Since
$\dim\,H=\infty$, this is not a compact operator. Thus, the assumed
property {\sf [H]} leads to a contradiction.
\end{example}

\vspace{.1in} Recall that, by \cite[Example 2]{b-g1}, $\dim_{\Bbb
K(H)}\Bbb K(H)=\dim\,H.$

Our following example shows that even a Hilbert $\Bbb K(H)$-module
${\mathscr M}$ such that $\dim_{\Bbb K(H)}{\mathscr M}<\infty$ need
not have property {\sf [H]}.

\begin{example}\label{primjer} (cf.~\cite[Example 1]{b-g1})
Let $H$ be a Hilbert space. For $\xi,\eta \in H$ define $\langle
\xi,\eta\rangle=e_{\xi,\eta} \in \Bbb K(H)$, where
$e_{\xi,\eta}(\nu)=(\nu|\eta)\xi$. Also, for $a \in \Bbb K(H)$,
define a left action on $\xi \in H$ in a natural way as the action
of the operator $a$ on the vector $\xi$.

In this way $H$ becomes a left Hilbert $\Bbb K(H)$-module. Notice
that the resulting norm coincides with the original norm on $H$.

We also know that $\dim_{\Bbb K(H)}H=1$. Indeed, if $\varepsilon$
is an arbitrary unit vector then each $\xi \in H$ admits a
representation of the form $\xi=\langle \xi, \varepsilon\rangle
\varepsilon$ (because $\langle \xi, \varepsilon \rangle
\varepsilon=e_{\xi,\varepsilon}(\varepsilon)=(\varepsilon|\varepsilon)\xi=\xi$).
This means that $\{\varepsilon\}$ is an orthonormal basis for $H$,
regarded as a $\Bbb K(H)$-module.

Notice that the entire preceding discussion was independent on the
(usual) dimension of the underlying space $H$. Suppose now that
$H$ is a separable infinite dimensional Hilbert space. We claim
that then $H$, as a Hilbert $\Bbb K(H)$-module, does not satisfy
{\sf [H]}.

To see this, let us fix an orthonormal basis $(\varepsilon_n)$ for
$H$. The sequence $(\varepsilon_n)$ is obviously norm-bounded.
Suppose that there exist a subsequence $(\varepsilon_{n_k})$ of
$(\varepsilon_n)$ and $\varepsilon_0 \in H$ such that
$\lim_{k\rightarrow \infty}\langle \varepsilon_{n_k},\xi \rangle
=\langle \varepsilon_0, \xi \rangle$ for all $\xi \in H$. In
particular, this would imply $\lim_{k\rightarrow \infty}\langle
\varepsilon_{n_k},\varepsilon_1 \rangle =\langle \varepsilon_0,
\varepsilon_1 \rangle$, i.e.,
$\|e_{\varepsilon_{n_k},\varepsilon_1}-e_{\varepsilon_0,
\varepsilon_1}\| \rightarrow 0$. But,
$\|e_{\varepsilon_{n_k},\varepsilon_1}-e_{\varepsilon_0,
\varepsilon_1}\|=\sup_{\|\eta\|=1}\|e_{\varepsilon_{n_k},\varepsilon_1}(\eta)-e_{\varepsilon_0,
\varepsilon_1}(\eta)\|
=\sup_{\|\eta\|=1}\|(\eta|\varepsilon_1)(\varepsilon_{n_k}-\varepsilon_0)\|=
\|\varepsilon_{n_k}-\varepsilon_0\|$ and the last expression
obviously does not converge to $0$ as $k\to \infty$.
\end{example}

\begin{remark}
Suppose that ${\mathscr M}$ is an arbitrary Hilbert $C^*$-module
over a $C^*$-algebra ${\mathscr A}$ of compact operators. It is well
known that there is a family $(H_j), j \in J$, of Hilbert spaces
such that ${\mathscr A}=\bigoplus_{j \in J}\Bbb K(H_j)$.
Furthermore, it then follows that ${\mathscr M}=\bigoplus_{j \in
J}{\mathscr M}_j$, where ${\mathscr M}_j=\overline{\Bbb
K(H_j){\mathscr M}}$ (i.e., ${\mathscr M}$ is an outer direct sum of
${\mathscr M}_j$'s, where each ${\mathscr M}_j$ is a full Hilbert
$\Bbb K(H_j)$-module).

Now, by \cite[Theorem 3]{b-g1} and the preceding example, we
conclude that if there exists $j_0 \in J$ such that
$\dim\,H_{j_0}=\infty$ then ${\mathscr M}_{j_0}$ cannot satisfy {\sf
[H]}. Consequently, ${\mathscr M}$ does not satisfy {\sf [H]}.
Namely, if $\dim\,{\mathscr M}_{j_0}=d$ (here $d$ can be an
arbitrary cardinal number), then, by Theorem 3 from \cite{b-g1},
${\mathscr M}_{j_0}$ is an orthogonal sum of $d$ copies of
$\mbox{}_{\Bbb K(H_{j_0})}H_{j_0}$, and, by Example \ref{primjer},
just one copy of $\mbox{}_{\Bbb K(H_{j_0})}H_{j_0}$ is enough to
ruin property {\sf [H]}.
\end{remark}

From the preceding discussion we conclude that if ${\mathscr M}$ is
a full Hilbert $C^*$-module satisfying {\sf [H]}, then ${\mathscr
M}$ is necessarily a Hilbert $C^*$-module over a $C^*$-algebra
${\mathscr A}$ of compact operators. Moreover, ${\mathscr A}$ has to
be of the form ${\mathscr A}=\bigoplus_{j \in J}\Bbb K(H_j)$ and
each $H_j$ must be finite dimensional. If, moreover, $J$ is of
finite cardinality, then ${\mathscr A}$ is finite dimensional. Next
we show that if $\mbox{card}\,J=\infty$ with $\dim\,H_j<\infty$ for
all $j \in J,$ then again ${\mathscr M}$ cannot satisfy {\sf [H]}.

First, in this situation, since $J$ as a set of an infinite
cardinality contains a countable subset $J^{\prime}$, ${\mathscr
M}=\bigoplus_{j \in J}{\mathscr M}_j$ can be written as the
orthogonal sum of the form ${\mathscr M}=\left(\bigoplus_{j \in
J^{\prime}}{\mathscr M}_j\right) \bigoplus \left(\bigoplus_{j \in
J\setminus J^{\prime}}{\mathscr M}_j\right)$. Thus, ${\mathscr M}$
contains, as an orthogonal summand, a submodule of the form
${\mathscr M}^{\prime}=\bigoplus_{n \in \Bbb N}{\mathscr M}_n$,
where each ${\mathscr M}_n$ is a module over $\Bbb K(H_n)$ and
$\dim\,H_n <\infty$. Moreover, each ${\mathscr M}_n$ is, by
\cite[Theorem 3]{b-g1}, unitarily equivalent to the orthogonal sum
of $d_n=\dim\,{\mathscr M}_n$ copies of $\mbox{}_{\Bbb K(H_n)}H_n,$
i.e., ${\mathscr M}_n\simeq\bigoplus_1^{d_n} \mbox{}_{\Bbb
K(H_n)}H_n$.

If we take just one copy of each $\mbox{}_{\Bbb K(H_n)}H_n$, we
conclude that ${\mathscr M}^{\prime}$ (and hence ${\mathscr M}$)
contains, as an orthogonal summand, a submodule of the form
${\mathscr M}^{\prime
\prime}\simeq\bigoplus_{n=1}^{\infty}\mbox{}_{\Bbb K(H_n)}H_n$. It
is now enough to prove that ${\mathscr M}^{\prime \prime}$ does not
satisfy {\sf [H]} and this can be argued essentially in the same way
as in Example~\ref{primjer}.

Observe that ${\mathscr M}^{\prime \prime}$ is also a Hilbert
$C^*$-module over a direct sum $\bigoplus_{n=1}^{\infty}\Bbb K(H_n)
\subset \Bbb K(H)$, where $H=\bigoplus_{n=1}^{\infty}H_n$ is an
infinite dimensional Hilbert space. For each $n \in \Bbb N$ take a
unit vector $\varepsilon_n \in H_n \subset H$. Let
$x_n=(\varepsilon_1,\varepsilon_2,\ldots,\varepsilon_n,0,0,\ldots
),\,n \in \Bbb N$. Notice that $\langle x_n,x_n
\rangle=\sum_{i=1}^ne_{\varepsilon_i,\varepsilon_i}$. Since this is
an orthogonal projection onto an $n$-dimensional subspace of $H$, we
have $\|x_n\|=1$; thus, $(x_n)$ is a norm-bounded sequence in
${\mathscr M}^{\prime \prime}$. Suppose now that there exists a
subsequence $(x_{n_k})$ and $x_0=(\xi_1,\xi_2,\ldots ) \in {\mathscr
M}^{\prime \prime}$ such that $\lim_{k \rightarrow \infty} \langle
x_{n_k},y\rangle = \langle x_0,y\rangle$ for all $y \in {\mathscr
M}^{\prime \prime}$. Inserting $y=(\varepsilon_1-\xi_1,0,0,\ldots )$
we obtain $\|\langle x_{n_k},y\rangle - \langle
x_0,y\rangle\|=\|\langle \varepsilon_1-\xi_1,\varepsilon_1-\xi_1
\rangle\| \rightarrow 0$, which implies $\xi_1=\varepsilon_1$.
Similarly, for $y=(0,\varepsilon_2-\xi_2,0,\ldots)$ we obtain
$\xi_2=\varepsilon_2$ and, proceeding in the same way,
$\xi_n=\varepsilon_n$ for all $n \in \Bbb N$. This gives us
$x_0=(\varepsilon_1,\varepsilon_2,\varepsilon_3,\ldots )$, which is
impossible since this sequence does not belong to ${\mathscr
M}^{\prime \prime}$.

\vspace{.1in}

After all, combining the preceding discussion with Proposition 2.1
from \cite{C-I-M-S}, we get our main result.

\begin{theorem}
A full Hilbert $C^*$-module over a $C^*$-algebra ${\mathscr A}$
satisfies condition {\sf [H]} if and only if ${\mathscr A}$ is a
finite dimensional $C^*$-algebra.
\end{theorem}

\vspace{.1in}

\begin{remark} We may ask ourselves if one could replace condition
{\sf[H]} with a weaker one:

\item{\sf[H']} For each norm-bounded sequence $(x_n)$ in
${\mathscr M}$ and for every $y\in {\mathscr M}$ there exists a
subsequence $(x_{n_k})$ of $(x_n)$ such that the sequence
$(\ip{x_{n_k}}{y})$ converges in norm.

\vspace{.1in}

Observe that {\sf[H']} is sufficient to prove Proposition~\ref{p1},
so full Hilbert $C^*$-modules with property {\sf[H']} have to be
over $C^*$-algebras of compact operators. Also, it is obvious that
{\sf[H']} is fulfilled in every Hilbert $C^*$-module over a finite
dimensional $C^*$-algebra. However, our next example shows that
{\sf[H']} does not characterize these Hilbert modules; in other
words, {\sf[H']} is not sufficient for {\sf[H]}.

Consider a separable infinite dimensional Hilbert space $H$ and the
$C^*$-algebra ${\mathscr A} \subset \Bbb K(H)$ of all diagonal (with
respect to a fixed orthonormal basis) operators with diagonal
entries converging to $0$. Let ${\mathscr M}={\mathscr A}$. Then
${\mathscr A}$ is a Hilbert $C^*$-module whose underlying
$C^*$-algebra ${\mathscr A}$ is infinite dimensional. By the
preceding theorem, the Hilbert $C^*$-module ${\mathscr A}$ cannot
satisfy {\sf[H]}.

On the other hand, since ${\mathscr A}$ is a Hilbert $C^*$-module
over the  (commutative) $C^*$-algebra ${\mathscr A}$ of compact
operators, by \cite[Theorem 4]{a} (see also its proof), all mappings
$T_y :{\mathscr A} \rightarrow {\mathscr A}$ given by
$T_y(x)=\langle y,x\rangle y$ are compact. But here we have
$T_y(x)=yx^*y=x^*y^2$ for all $y\in {\mathscr A}.$ In particular,
taking self-adjoint $y$ we get that $x\mapsto x^*y$ is compact for
every positive $y\in {\mathscr A},$ and since positive elements of a
$C^*$-algebra span the whole $C^*$-algebra, we get that the operator
$x\mapsto x^*y=\ip{y}{x}$ is compact for every $y\in {\mathscr A}.$
This shows that our Hilbert $C^*$-module ${\mathscr A}$ satisfies
{\sf[H']}.
\end{remark}

\vspace{.1in}

{\it Acknowledgement.} The authors are thankful to the referee for
several valuable suggestions for improving and clarifying the
original manuscript.

\end{document}